\documentclass[12pt,twoside,reqno]{amsart}
\usepackage[colorlinks=true,citecolor=blue]{hyperref}
\usepackage{mathptmx, amsmath, amssymb, amsfonts, amsthm, mathptmx, enumerate, color,mathrsfs}
\usepackage{graphicx}

\usepackage{multirow}
\usepackage{epstopdf}
\usepackage{multicol}
\usepackage{algorithm}
\usepackage{algorithmic}
\usepackage{epstopdf}

\newtheorem{theorem}{Theorem}[section]

\theoremstyle{definition}

\newtheorem{remark}[theorem]{Remark}
\numberwithin{equation}{section}

\normalsize

\newcommand{\Z}{Z\!\!\!Z}
\newcommand{\C}{I\!\!\!\!C}

\newcommand\meas{\mathop{\rm meas}}

\def\C{{\mathcal C}} 

 \def\al{{\alpha}}

\newcommand{\q}{\quad} \newcommand{\qq}{\qquad}

\newtheorem{Theorem}{Theorem}[section]
\newtheorem{Proposition}{Proposition}[section]
\newtheorem{Remark}{Remark}[section]

\newtheorem{Definition}{Definition}[section]
\newtheorem{Example}{Example}[section] \newtheorem{Assumption}{Assumption}[section]

\newcommand{\bd}{\begin{displaymath}}
\newcommand{\ed}{\end{displaymath}}
\newcommand{\be}{\begin{equation}}
\newcommand{\ee}{\end{equation}}
\newcommand{\bea}{\begin{eqnarray}}
\newcommand{\eea}{\end{eqnarray}}
\newcommand{\bda}{\begin{eqnarray*}}
	\newcommand{\eda}{\end{eqnarray*}}
\newcommand{\ba}{\begin{array}}
	\newcommand{\ea}{\end{array}}

\newcommand{\B}{I\kern -.35em B}

\newcommand{\R}{I\kern -.35em R}

\newcommand{\To}{\rightrightarrows}

\newcommand{\st}{\subset}

\newcommand{\LL}{{\mathcal L}}

\newcommand{\W}{{\mathcal W}}

\newcommand{\e}{\varepsilon}

\newcommand{\sth}{ \, :\;}

\newcommand{\dd}{\mbox{\rm\,d}}

\newcommand{\X}{{\mathcal X}}

\renewcommand{\Z}{{\mathcal Z}}
\newcommand{\F}{{\mathcal F}}
\renewcommand{\O}{{\mathcal O}}

\def\R{\mathbb{R}}
\newcommand{\reff}{\eqref}
\newcommand{\bino}{\bigskip\noindent}

\newcommand{\Dex}{{\Delta x}}

\def\bb{\textcolor{blue}}

\begin{document}
\setcounter{page}{1}

\vspace*{2.0cm}
\title[Strong Metric
Subregularity]
{Strong Metric
	Subregularity of the Optimality Mapping and  Second-Order Sufficient Optimality Conditions in Extremal
	Problems with Constraints}

\begin{center}
	{\it Dedicated to Vladimir M. Tikhomirov}	
	
\end{center}

\author[N. Osmolovskii, V. Veliov]{ Nikolai P. Osmolovskii$^{1,*}$,  Vladimir M. Veliov$^2$}
\maketitle
\vspace*{-0.6cm}

\begin{center}
{\footnotesize

$^1$Systems Research Institute, Polish Academy of Sciences,  Warsaw, Poland\\
$^2$Institute of Statistics and Mathematical Methods in Economics,
TU Wien, Austria
}\end{center}

\vskip 4mm {\footnotesize \noindent {\bf Abstract.} This is a review paper, summarizing without proofs recent results by the authors 
	on the property of strong metric subregularity (SMSR) in optimization.
	It presents sufficient conditions for SMSR of the optimality mapping 
	associated with a set of necessary optimality conditions in three types of constrained optimization problems: 
	mathematical programming, calculus of variations, and optimal control. The conditions are based on 
	second-order sufficient optimality conditions in the corresponding optimization problems and guarantee 
	small changes in the optimal solution and Lagrange multipliers for small changes in the data.


  }

\renewcommand{\thefootnote}{}
\footnotetext{ $^*$Corresponding author.
\par
E-mail address: nikolai.osmolovskii@ibspan.waw.pl (N. Osmolovskii), vladimir.veliov@tuwien.ac.at (V. Veliov).
\par
Received xx, x, xxxx; Accepted xx, x, xxxx.
\rightline {\tiny   \copyright  2022 Communications in Optimization Theory}}

\section{Introduction}

This paper reviews results on Lipschitz stability under perturbations in extremal 
problems with constraints. The analysis is based on (and restricted to) the abstract concept of 
{\em strong metric subregularity} (SMSR) for set-valued mappings in Banach spaces. 
This concept can be formalized
in several ways depending on the topologies used in spaces.
To set the stage, let $(\X,\| \cdot \|_\X )$ be a Banach space, in which one more norm
$\| \cdot \|_\X'$ is defined such that:
$$
\| s \|'_\X \leq \| s \|_\X\q  \forall \,s \in \X,
$$ 
that is, $\| \cdot \|'_\X$ is a weaker norm than the basic 
norm $\| \cdot \|_\X$.
For example,  
$$
\X=L^\infty(0,1),\q \| \cdot \|_\X=\| \cdot \|_\infty, \q
\| \cdot \|_\X'=\| \cdot \|_2.
$$
Let $\LL : \X \To \Z$
be a set-valued mapping, where $(\Z,\| \cdot \|_\Z )$ is a normed  space. 

\begin{Definition}\label{def1} 
	The mapping $\LL : \X \To \Z$ has the property SMSR at {$\hat s$ for $\hat z$} if $\hat z\in \LL(\hat s) $ 
	and there exist neighborhoods $\O_{\hat s} \ni \hat s$ (in the norm $\| \cdot \|_\X$), 
	$\O_{\hat z} \ni \hat z$ and a number $\kappa$ such that the relations 
	\bd
	s \in \O_{\hat s}, \quad z \in \O_{\hat z}, \quad z \in \LL(s)
	\ed
	imply that 
	\be \label{Esmsr}
	\| s - \hat s \|'_\X \leq \kappa \| z - \hat z \|_{\Z}.
	\ee
\end{Definition}

We will use this definition for $\hat z=0$.

The definition of SMSR that we use is a slight extension 
of the standard one, 
introduced under this name in \cite{AD+TR-04}; see also \cite[Chapter~3.9]{AD+TR-book2}
and the recent paper \cite{Cibulka+Dontchev+Kruger-18}. 
Notice that in the above definition, the smaller norm, $\| \cdot \|'_\X$,  is involved in \reff{Esmsr}, 
while the neighborhood $\O_{\hat s}$ is with respect to the larger one. The difference with the standard
definition in \cite{AD+TR-04} and \cite{Cibulka+Dontchev+Kruger-18} is that the two norms in $\X$
coincide there. We mention that in some cases (including the last section of the present paper)
it makes sense to involve two norms also in the space $\Z$\bb{ (see \cite{MQ+VV:SICON-13,Corella+Quinc+Vel-20})}.

Other versions of the SMSR property were introduced and utilized in 
\cite{Bonnans-94,Bonn+Shap-book-2000,Klatte+Kummer-2002-book}. 
It is well recognized that the SMSR property of the mapping associated with the first order 
optimality conditions of optimization problems, the so-called {\em optimality mapping}, 
is a key property for ensuring convergence with error estimates 
of numerous methods for solving optimal control problems: discretization methods, gradient methods,  
Newton-type methods, etc. (see e.g.
\cite{Bonnans-94,Alt+S+S-16,Pre_Sca_Vel_MReg-17,Cibulka+Dontchev+Kruger-18,Ang+Cor+Vel_MPS_21}, 
in addition to a large number of papers where the SMSR property is implicitly used).

The results presented in the paper focus on the SMSR property of the system of necessary optimality 
conditions of a general mathematical programming problem in Banach spaces and several classes
of optimal control problems for ordinary differential equations. The related issue of sufficient optimality
conditions is also addressed. For results on the SMSR property of the optimality mapping of 
PDE-constrained optimization problems we refer to the recent papers 
\cite{DJV-22-elliptic,Corella+Jork+Vel_COAP_23}.     

The paper is organized as follows. In Section \ref{sec2} we discuss a  mathematical programming (MP) 
problem in a two-norm space,  define the  Mangasarian-Fromovitz constraint qualification condition, 
introduce assumption about the ``two-norm differentiability'' of data, recall the  Karush-Kuhn-Tucker (KKT)  
optimality conditions and sufficient conditions of the second order  in the MP problem, and, finally, 
we formulate a sufficient condition for the  SMSR property for the KKT mapping in the MP problem. 
The last mapping is related to the KKT conditions. In Section \ref{sec3} we apply the abstract results 
of Section \ref{sec2} to a general calculus of variations problem represented as a Mayer-type optimal 
control problem without the control constraints. In Section \ref{sec4} we discuss the second-order 
sufficient condition for a weak local minimum and the SMSR property in a Mayer-type optimal control 
problem with control constraints $u\in U$, but without endpoint constraints. 
We assume that the set $U$ is defined 
by a finite number of smooth inequalities with linear independent gradients of the active constraints.

\section{ Mathematical programming problem in a space with two norms}\label{sec2}

Let $X$ and $Y$ be Banach spaces. Consider the problem
$$ 
\min\, \varphi(x), \q  f(x) \le 0,\q \q g(x)=0, \eqno{\rm (MP)}
$$
where $\varphi:X \to \R$, $f=(f_1,\ldots,f_m)^*:X \to \R^m$, $g : X \to Y$ are $C^1$ mappings. 
The symbol $*$ denotes transposition; when applied to a space it denotes the dual space.

Let $\hat x\in X$ satisfy the constraints
$f( x) \le 0$, $g( x)=0$. Denote by 
$$ I=\{i \sth f_i(\hat x)=0 \}$$
the set of active indices.
The following assumption is called Mangasarian-Fromovitz constraint qualification$^3$ \footnote{%
$^3${Usually the Mangasayan–Fromowitz condition is equivalently formulated as follows: $g'(\hat x)$ is surjective,
and there exists $\xi \in X$ such that $g'(\hat x) \xi = 0$ and $f_i'(\hat x)\xi < 0$, $i \in I$.}
}
 (MFCQ).

\begin{Assumption}\label{assump1}    $g'(\hat x)X=Y$ and the gradients $f'_i(\hat x)$, $i\in I$, are  
	positively independent on the subspace $\ker g'(\hat x)$.
\end{Assumption} 

Recall that the  functionals $l_1,\ldots ,l_k\in X^*$ are  {\it positively independent on a subspace 
	$L\subset X$} if 
$$\sum_{i=1}^k \lambda _i l_i(x)=0 \q \forall\, x\in L,\q \lambda_i\ge0\q \forall\,i\q  
\Longrightarrow \q  \lambda_i=0 \q \forall \,i.$$ 

Let $\R^{m*}$ denote the space of row vectors of dimension $m$.

\begin{Theorem}\label{th2.1} If $\hat x$ is a local minimizer in problem (MP) and Assumption \ref{assump1} holds, 
	then  there exist $\hat\lambda\in \R^{m*}$ and $\hat y^*\in Y^*$ such that 
	\bda
	&&  \varphi'(\hat x) +  \hat\lambda f'(\hat x)+ \hat y^*g'(\hat x)  = 0, \\ 
	&&  \hat\lambda f(\hat x)  = 0,  \q  \hat\lambda\ge 0, \\ 
	&&  g(\hat x) = 0, \q  f(\hat x) \le 0. 
	\eda
\end{Theorem}

The above relations are known as Karush-Kuhn-Tucker (KKT) system, and the triple 
$(\hat x,\hat\lambda,\hat y^*)$ that satisfies these conditions is called the KKT point.

In what follows we fix the  KKT point $(\hat x,  \hat \lambda, \hat y^*)$. 

\begin{remark} Theorem 2.1 is actually the rule of Lagrange multipliers in the normal case (where $\lambda_0 = 1$). 
{It is attributed to Karush-Kuhn-Tucker, although they have not proved the theorem in the above form.
Their proof involves an {\it assumption}, 
which is not easier to verify than solving the original problem itself.} (This is precisely why Mangasarian and Fromowitz 
introduced their condition!) The history here is long, dating back to the 18th century, but eventually the proof was 
given by Mangasarian and Fromowitz in 1967; it easily follows also from the  earlier work of
Dubovitskii and Milyutin in 1965 \cite{DM-1965}.
\end{remark}

Assume that in the Banach space $(X,\|\cdot\|)$ there is another, weaker norm $\|\cdot\|'$, that is 
$\|x\|'\le \|x\|$ for all $ x\in X$.

We make the following ``two-norm differentiability'' assumptions for the mappings $\varphi$, $f$ and $g$. 
First we formulate it for $g$.

\begin{Assumption}\label{assump2}
	There exists a neighborhood $\hat O$ of $\hat x$ (in the norm $\| \cdot \|$) such that the following 
	conditions are fulfilled for all $\Delta x\in X$ such that $\hat x + \Delta x \in \hat O$:
	\begin{itemize}
		
		\item[(i)] The operator $g$ is continuously Fr\'echet differentiable in $\hat O$ in the norm 
		$\| \cdot \|$, and the derivative $g'(\hat x)$ is continuous operator 
		with respect to the norm $\| \cdot \|'$; moreover, the following representation holds true 
		$$
		g(\hat x + \Delta x) = g(\hat x) + g'(\hat x) \Delta x + r(\Dex), 
		$$
		with   
		$$
		\|  r(\Dex) \|_Y \leq \theta(\| \Dex \|) \| \Dex \|',$$ 
		where $\theta(t)\to0$ as $t\to0+$. 	
		
		\item[(ii)]	  There exists a bilinear mapping $Q: X \times X \to Y$ 
		such that   
		$$
		g'(\hat x + \Dex)= g'(\hat x) + Q( \Dex,\cdot) + \bar r(\Dex), 
		$$
		$$
		\|Q(x_1,x_2)\|_Y \leq C \| x_1 \|' \| x_2 \|',
		$$ 
		where $C$ is a constant and $\bar r$ satisfies 
		$$
		\sup_{\| x \|' \leq 1} \| \bar r(\Dex)(x) \|_Y        \leq \theta(\| \Dex \|) \| \Dex \|'. 
		$$
		
		\item[(iii)] $\varphi$ and $f$  satisfy similar conditions with bilinear mappings $Q_0$ and $Q_f$.
	\end{itemize}
\end{Assumption}

It is easy to show that (i) and (ii)  imply the expansion
$$
g(\hat x+\Dex) = g(\hat x) + g'(\hat x)\Dex +\frac{1}{2}Q(\Dex,\Dex) + \hat r(\Dex),
$$
where 	 
$$
\|  \hat r(\Dex) \|_Y \leq \theta(\| \Dex \|) (\| \Dex \|')^2.  
$$
Similarly,  $\varphi$ and $f$ have the following expansions
$$
\varphi(\hat x+\Dex) = \varphi(\hat x) + \varphi'(\hat x)\Dex +\frac{1}{2}Q_0(\Dex,\Dex) + \hat r_0(\Dex),
$$
$$
f(\hat x+\Dex) = f(\hat x) + f'(\hat x)\Dex +\frac{1}{2}Q_f(\Dex,\Dex) + \hat r_f(\Dex),
$$
where $\hat r_0$, $ \hat r_f$ satisfy properties similar to that of $ \hat r$.

Define the quadratic functional $\Omega:X \to \R$ as
$$
\Omega(x)  :=  Q_0(x,x) + \hat \lambda Q_f(x, x) + \hat y^*\, Q(x,x),
$$
and the so-called {\em critical cone} at the point $\hat x$  as
\bd
K :=  \{ x \in X \sth \varphi'(\hat x) x\le 0, \q  f'_i(\hat x) x \le 0\; \mbox{ for }\; i \in I,\q  g'(\hat x) x = 0\}.
\ed

\begin{Assumption}\label{assump3}
	There exists a constant $c_ 0 > 0$ such that 
	\bd
	\Omega(x) \geq c_0\, (\| x \|')^2 \quad \forall \, x \in K.
	\ed
\end{Assumption}

\begin{Theorem} 
	Let Assumptions \ref{assump1}--\ref{assump3} be fulfilled.
	Then the following quadratic growth condition for $\varphi$ holds:
	there exist $c>0$ and $\varepsilon>0$ such that  
	\bd
	\varphi(x )- 	\varphi(\hat x) \ge c \, (\| x-\hat x\|')^2
	\ed 
	for all admissible $x$ satisfying $\|x-\hat x\| < \varepsilon$. Hence,  $\hat x$ is a strict 
	local minimizer in the problem. 
\end{Theorem}

For the  KKT point $(\hat x,  \hat \lambda, \hat y^*)$ 
we	split the set of active indices $I$ into two parts:
$$
I_0=\{ i\in I : \, \hat \lambda_i=0\}, \qq I_1 = \{i \in I\sth   \hat \lambda_i >0 \}.
$$

\begin{Assumption}\label{assump4}(strict MFCQ) 
	The image  $g'(\hat x) X$ is closed and 
	$$\ba{rcl}&&
	\lambda_i \ge 0  \q \forall\, i \in I_0, \q \sum_{i\in I}  \lambda_i  f_i'(\hat x)+y^* g'(\hat x) = 0 \\\\&& \Longrightarrow 
	\;\lambda_i = 0\q  \forall\,i \in I, \q y^* = 0.\ea$$
\end{Assumption}

Strict MFCQ (Assumptions \ref{assump4}) implies   MFCQ (Assumptions \ref{assump1}). 
In particular, {it} implies the condition $g'(\hat x) X=Y$. 
Moreover,  strict MFCQ  implies that 
$(\hat x,  \hat \lambda,\hat y^*)$ is the unique KKT point with the fixed $\hat x$.

The KKT system can be formulated in a more compact way as 
$$
L'_x(x,\lambda,y^*) = 0, \q
g(x) = 0, \q
f(x) \in N_{\R^m_+}(\lambda),
$$
where $L(x,\lambda,y^*) =  \varphi(x)  +   \lambda  f(x) + y^* g(x) $  is the Lagrangian, and 
\bd
N_{\R^m_+}(\lambda) := \left\{ \begin{array}{cl} 
	\big\{\alpha \in \R^m  \sth  \langle \alpha, \beta -\lambda \rangle \leq 0
	{\text{ for all } \beta \in \R^m_+}\,\big\} & \text{if } \lambda \in \R^m_+,\\ \\
	\emptyset &  \text{if } \lambda \not\in \R^m_+ \end{array}\right.
\ed
is the normal cone to $\R^m_+$ at $\lambda \in \R^m$. Obviously,  the condition  $  f(x) \in N_{R^m_+}(\lambda)$ 
incorporates the {\em complementary slackness} condition $\lambda f(x) = 0$.

Define the KKT mapping $\F: X\times \R^{m*}\times Y^*\To \R^m \times Y\times X^* $:
$$
\F(x, \lambda, y^*) := \left( \begin{array}{c}
f(x) - N_{\R^m_+}(\lambda)\\ g(x) \\
L'_x(x, \lambda, y^*) 
\end{array} \right).
$$

The elements of the image space $\R^m \times Y\times X^*$ are  denoted by $(\xi,\eta,\zeta)$. 	
Then the KKT system is equivalent to $$(0,0,0)\in \F(x, \lambda, y^*).$$

Further, set
$$X'' = \{x^* \in X^*  \sth x^* \mbox{ is continuous with respect to } \| \cdot \|' \}.$$

For any  $x^*\in X''$, we set
\bd
\| x^* \|'' :=\sup_{\|x\|' \leq 1}  | x^*(x)  |.
\ed
Obviously, $\| x^* \|'' \ge \| x^* \|$ for every $x^* \in X''$.

Define the space $$\X := X   \times \R^{m*}\times Y^*$$ with the following two norms: for 
$	s = (x,\lambda, y^*) \in \X$
$$
\| s \|_\X := \| x \| +|\lambda|+ \| y^* \| \quad \mbox{ and } \;\;  \| s \|'_\X :=  \| x \|'+|\lambda|+ \| y^* \| .$$
Also define {$\Z := \R^m\times Y \times X''$} 
with the following norm: for $z = (\xi, \eta, \zeta) \in \Z$
\bd
\| z \|_\Z =  | \xi |  + \| \eta \| + \| \zeta \|''.
\ed
Observe that due to Assumption 2  we have $f'_i(\hat x) \in X''$ and $ y^* g'(\hat x) \in X''.$ 	 
Thus $\F: \X \To \Z$. 

\begin{Theorem}\label{th3}
	Let $\hat s = (\hat x,  \hat \lambda,\hat y^*)\in\X$ be a KKT point for problem (MP), and let assumptions 
	\ref{assump2}--\ref{assump4} 
	be fulfilled at this point. Then the KKT mapping $\F$ 
	has the property SMSR at $\hat s$ for zero. More precisely, there exist  constants $a > 0$, $b > 0$ 
	and $\kappa \geq 0$ 
	such that for every $z = (\xi,  \eta, \zeta )\in\Z$  such that 
	$ 	\|z\|_\Z \le b$ and for every solution $s = (x,\lambda,y^*)\in\X$ of the inclusion $z \in \F(s)$
	with $\| x - \hat x\| \leq a$ it holds that 
	$$ \|s-\hat s\|_\Z'=
	\| x - \hat x \|' +| \lambda  - \hat \lambda |+ \| y^* - \hat y^* \|  
$$
$$  \leq 
	{\kappa (| \xi |  + \| \eta \| +  \|\zeta \|'' )}=\kappa\|z\|_\Z.
$$        
\end{Theorem}

\bino
Proofs of the results in this section are given in  \cite{Osmol+Vel-21-JMAA} and \cite{Osmol+Vel-21-C-and-C}.

\section{General problem of  calculus of variations}\label{sec3}

The obtained abstract subregularity result can be applied to the following 
optimal control problem of Mayer's type on a fixed time interval $[t_0,t_1]$:
$$
\mbox{minimize}\q \varphi_0(x(t_0),x(t_1))
$$
subject to
$$ 
\dot x(t) = f(x(t),u(t)), \qquad t \in  [t_0,t_1], 
$$
$$ 
\psi_j(x(t_0),x(t_1))=0, \q j=1,\ldots,s,
$$
$$  
\varphi_i(x(t_0),x(t_1))\le0, \q i=1,\ldots,k,
$$
where  $f:\R^n\times\R^m\to\R^n$, $\psi_j:\R^n\times\R^n\to\R$, 
$j=1,\ldots,s$, $\varphi_i:\R^n\times\R^n\to\R$, $i=0,\ldots,k$, are $C^2$ mappings. We call this 
briefly the {\it Mayer problem}.

We stress that control constraint  $u\in U$ is not involved, therefore, 
we actually deal with a {\em general problem of calculus of variations}. 
The strict  Mangasarian-Fromovitz condition and the coercitivity condition take specific forms, 
in which two norms are used similarly as, for example, in the papers by Malanowski and Maurer 
\cite{Mal-94},\cite{Mal-Ma-96},\cite{Mal-01}.  
The main novelty of the strong subregularity result is that the coercivity condition is posed on a 
{\em critical cone} {with} {initial/terminal constraints}  of equality and inequality type 
for the state variable.

We consider the Mayer problem for trajectory-control pairs  
$w(\cdot)=(x(\cdot), u(\cdot))$ with  measurable and essentially bounded $u:[t_0,t_1]\to\R^m$ and 
absolutely continuous $x:[t_0,t_1]\to\R^n$. Thus the admissible points in the problem belong to the space
\bd
\W := W^{1,1}([t_0,t_1];\R^n)\times L^\infty([t_0,t_1];\R^m) 
\ed
with the norm
\bd
\|w\|=\|x\|_{1,1}+\|u\|_\infty.
\ed
We introduce a weaker norm
$$ \|w\|'=\|x\|_{\infty}+\|u\|_2.$$
Recall that 
$   \|u\|_2\le c \|u\|_\infty$ and $ \|x\|_{\infty}\le c\|x\|_{1,1}$ with some $c>0$, 
therefore the norm $\|w\|'$ is really weaker than $\|w\|$. Of course, we could preserve the norm $ \|x\|_{1,1}$ 
in the definition of  $\|w\|'$, but for the coercitivity Assumption \ref{assum3.2} (see below) and some estimates 
this weaker norm for $x$  is more convenient.

Set
$$ q=(x(t_0),x(t_1)),\q \alpha= ( \alpha_1,\ldots, \alpha_k)\in\R^{k*}.$$
We introduce the  {\it {Hamiltonian} (Pontryagin function) } and the {\it endpoint Lagrange function}:
$$
H(x,u,p) = p \,f(x,u), \q  l(q,\alpha_0,\alpha,\beta)= \sum_{i =0}^k \alpha_i \varphi_i(q) + \beta \psi(q),
$$
where $p\in\R^{n*}$, $\alpha\in \R^{k*}$,  $\beta\in \R^{s*}$ are row vectors, $\alpha_0$ is a number.
The local form of the Pontryagin maximum principle for the problem reads as follows. 

\begin{Theorem}\label{th3.1} If $ w=( x, u) \in \W$ is a local minimizer in the Mayer problem, then 
	there exists a nontrivial tuple $(\alpha_0,\alpha,\beta, p)\in \R\times \R^{k*}\times \R^{s*} 
	\times  W^{1,1}([t_0,t_1];\R^{n*})$	such that 
	$$
	\alpha_0\ge0,\q    \alpha\ge0, \q   \alpha\varphi(q)=0,
	$$
	$$  
	- \dot p =  H_x( w,p), \q H_u( w,p) = 0, \q (-p(t_0), p(t_1))  =  l'_q(q,\alpha_0,\alpha,\beta).
	$$
\end{Theorem}

Theorem \ref{th3.1}  could also be called the {\em general Euler–Lagrange equation}. Just as Theorem \ref{th2.1} it also follows from the 1965  work \cite{DM-1965} of Dubovitskii and Milyutin .

Let $\hat w=(\hat x, \hat u)$ be an admissible point in the Mayer problem and the tuple
$(\hat w, \hat\alpha_0,\hat\alpha,\hat\beta, \hat p)$ satisfies the conditions in the theorem.
{Denote
	$$
	I=\{i\in\{1,\ldots,k\}\sth \varphi_i(\hat q)=0  \}, \q I_0=\{i\in I  \sth \hat \alpha_i=0\}.
	$$
}
The next assumption is an analogue of the strict Mangasarian--Fromovitz condition.

\begin{Assumption}\label{assum3.1} 	The relations
	$$
	\alpha_0=0,\q (\alpha,\beta, p)\in \R^{k*}\times \R^{s*} \times W^{1,1} ([t_0,t_1];\R^{n*})  ,\q      
	\alpha_i \ge 0 \q \mbox{ for all} \q i  \in I_0,   
	$$ 
	$$ 
	- \dot p = H_x(\hat w,p), \q H_u(\hat w,p) = 0,  \;\,\; (-p(t_0), p(t_1))  =  l'_q(\hat q,0,\alpha,\beta) 
	$$
	imply that $\al = 0$, $ \beta = 0$,  $p = 0$.
\end{Assumption}

{Assumption \ref{assum3.1} implies that, for a given $\hat w=(\hat x, \hat u)$, the tuple 
	$(\hat \alpha_0,\hat\alpha,\hat\beta, \hat p)$ satisfying these first order necessary optimality conditions  
	is unique and $\hat \alpha_0 > 0$. Therefore, one can always set $\hat \alpha_0 = 1$.
}
For brevity we denote
{$$
	l_{q} (\hat q, 1, \hat \al,\hat \beta)=l_{q} (\hat q),\q	
	l_{qq} (\hat q, 1, \hat \al,\hat \beta)=l_{qq} (\hat q).	
	$$}	
For this tuple, we define the quadratic form
\bd
\Omega(w) = \langle l_{qq} (\hat q) q, q \rangle
+ \int_{t_0}^{t_1} \langle H_{ww}(\hat w,\hat p)(t) w(t), w(t) \rangle  \dd t,
\ed
where $w(t)=(x(t),u(t))$, $q=(x(t_0),x(t_1))$, and
\bd \ba{rcl}&& 
\langle H_{ww}(\hat w,\hat p)\, 
w, w \rangle\\& =&\langle H_{xx}(\hat w,\hat p) x, x \rangle  +
2 \langle H_{ux}(\hat w,\hat p) x, u \rangle +\langle H_{uu}(\hat w,\hat p)  u, u \rangle.  \ea
\ed

For the pair  $\hat w=(\hat x, \hat u)$,  define the critical cone as the set of all those 
pairs $( x, u)\in \W$ which satisfy 
$$ 
\dot x = f_x(\hat w) w,  \q \psi'(\hat q) q=0, \q 
\varphi'_i(\hat q)  q \le0,\q i\in I\cup\{0\}.
$$

\begin{Assumption}\label{assum3.2}
	There exists a constant $c_ 0 > 0$ such that 
	\bd
	\Omega(w) \geq c_0\, (\| w \|')^2 \quad \forall \, w \in K.
	\ed
\end{Assumption}

As in the case of mathematical programming problem, this assumption implies a weak local minimum 
at the point $\hat w$, and even more, it implies a quadratic growth condition for the cost in a neighborhood 
of the point $\hat w$: $$\varphi_0(q)-\varphi_0(\hat q) \ge c (\| w-\hat w \|')^2
$$
with some $c>0$. 

Note that the assumptions about `` differentiability with respect to two norms '' of problem mappings made 
in the abstract problem (MP) are automatically fulfilled.

Introduce the space 
$$
\X=\big\{s=(x,u,p,\al,\beta )\in W^{1,1}\times L^\infty\times  W^{1,1}\times \R^{k*}\times \R^{s*}\big\}.
$$
Consider the space  of disturbances
\bd
\Z = \big\{z= (\pi,\rho, \nu,  \eta, \mu, \xi) \in L^1 \times L^2 \times \R^{2n*} \times  L^1 \times \R^{s} \times 
\R^{k} \big\}.
\ed	
and a perturbed system of optimality conditions
$$
\dot p + p \,f_x(w) =\pi,$$ 
$$               p f_u(w) = \rho,$$ 
$$        (p(t_0), -p(t_1)) +  \varphi'_{0}(q) + \al \varphi'(q)  + \beta \psi'(q)=\nu,$$
$$f(x,u)-\dot x=\eta,$$
$$\psi(q)=\mu,$$       
$$\varphi(q)\le\xi, \q \alpha(\varphi(q)-\xi)=0.$$

Theorem \ref{th3}, obtained for  abstract problem (MP), implies the following result.

\begin{Theorem}\label{th3.2} 
	Let $\hat s = (\hat x, \hat u, \hat p, \hat\alpha, \hat \beta) \in \X$ be a solution 
	of the original optimality system 
	and let assumptions \ref{assum3.1}-\ref{assum3.2}   be fulfilled. 
	Then there exist  constants $a > 0$, $b > 0$ and $\kappa \geq 0$ 
	such that for every $z = (\pi,\rho,\nu, \eta, \mu, \xi) \in  \Z$ with 
	$$\| z \|_\Z := \| \pi \|_1  + \| \rho \|_2+ | \nu | + \| \eta \|_1 + | \mu | + | \xi | \leq b$$ 
	and for every solution $s = (x,u , p, \al, \beta) \in \X$ of the perturbed system
	with $$\| x - \hat x\|_{1,1} + \| u - \hat u \|_\infty \leq a\,$$
	 it holds that 
	$$
	\| x - \hat x \|_{1,1} + \| u - \hat u \|_2 +  \| p - \hat p \|_{1,1} + | \al  - \hat \al |+ | \beta  - \hat \beta |   
	\leq \kappa \| z \|_\Z.
	$$
\end{Theorem}

Proofs of the results in this section are given in  \cite{Osmol+Vel-21-JMAA}.	


\section{The simplest optimal control problem} \label{sec4} 

Here we investigate the following Mayer-type optimal 
control problem  without endpoint constraints but with control constraint $u\in U$:

\be\label{1}
\mbox{minimize} \q J(x,u):=F(x(0),x(1)), 
\ee 
\be\label{2} 
\dot x(t)=f(x(t),u(t))\q \mbox{a.e. in}\q [0,1], 
\ee
\be\label{3} 
G(u(t))\le0 \q \mbox{a.e. in}\q [0,1], 
\ee
where $F:\R^{2n}\to\R$, $f: \R^{n+m}\to\R^n$, and $G:\R^m\to\R^k$  are of class 
$C^2$,  $u\in L^\infty$, $x\in W^{1,1}$.

Again, we investigate the property of {\em strong metric subregularity}  (SMSR) of the 
optimality mapping,
associated with the system of first order necessary optimality conditions (Pontryagin's conditions in local form) 
for problem (\ref{1})--(\ref{3}).

{According to (\ref{3}), the set of admissible control values is}
$$ 
U:=\{v\in \R^m: \; G(v)\le0\}.
$$
Let $G_i$ denote the $i$th component of the vector $G$. For any $v\in U$ define the set of active indices
$$    
I(v)=\{i\in\{1,\ldots,k\}\,:  \q G_i(v)=0\}.
$$

\begin{Assumption}\label{assum4.1}{\em (regularity of the control constraints)} The set $U$ is nonempty and	
	at each point $v\in U$ the gradients $G_i'(v)$, $i\in I(v)$, are linearly independent.
\end{Assumption}

As in the previous section, we use the notations 
$$ 
q=(x(0),x(1))=(x_0,x_1),\q w=(x,u), \q \W=W^{1,1}\times L^\infty.
$$

Let $\hat w=(\hat x,\hat u)\in\W$ be the point examined for optimality, and let $\hat q=(\hat x(0),\hat x(1))$ . 
By $\hat \lambda\in L^\infty$ we denote the Lagrange multiplier, which relates to the control constraint (\ref{3}), 
and by $\hat p\in W^{1,1}$ we denote the adjoint variable, which relates to the control system (\ref{2}).

\begin{Assumption}\label{assum4.2}
	{The triplet $(\hat w,\hat p, \hat \lambda) \in \W \times  W^{1,1} \times L^\infty$ satisfies the following 
		system of equations and inequalities:}

\be\label{3a} \hat \lambda(t)\ge0, \q \hat \lambda(t)  G(\hat u(t))=0\q \mbox{a.e. in}\q [0,1],\ee
\be\label{4}( -\hat p(0), \hat p(1))=F'(\hat q ), \ee
\be\label{5} \dot {\hat p}(t) +\hat p(t)\, f_x(\hat w(t))=0 \q \mbox{a.e. in}\q [0,1],\ee
\be\label{6} \hat p(t)\,f_u(\hat w(t))+\hat\lambda(t) G'(\hat u(t))=0 \q \mbox{a.e. in}\q [0,1],\ee
\be\label{9} -\dot{\hat x}(t)+f(\hat w(t))=0 \q \mbox{a.e. in}\q [0,1],\ee
\be\label{10}  G(\hat u(t))\le0 \q \mbox{a.e. in}\q [0,1].\ee
\end{Assumption}

Observe that this system represents the  first order necessary optimality condition for a weak local minimum
of the pair $\hat w=(\hat x,\hat u)$, which is a local minimum in the space $\W$. 	The latter means that  
$J(\hat x,\hat u) \leq J(x,u)$ for every admissible pair $(x,u)$ which is close enough to $(\hat x, \hat u)$
in the space $\W$. Later on, we refer to the system \reff{3a}--\reff{10} as to {\em optimality system}. 
Namely, if $\hat w$ is a point of weak local minimum in problem (\ref{1})--(\ref{3}),
then there exist $\hat p\in W^{1,1}$ and $\hat \lambda \in L^\infty$  such that the optimality system is fulfilled.
Note that for a given $\hat w$ the pair $(\hat p,\hat\lambda)$ is uniquely determined by these conditions.

Introduce the  {\em Pontryagin function} (often called the Hamiltonian) and the {\em augmented Pontryagin function}
$$ 
H(w,p)=p\, f(w), \q \bar H(w,p,\lambda)=p\, f(w)+\lambda\, G(u).
$$
Then equations (\ref{5}) and  (\ref{6}) take the form
$$ 
-\dot {\hat p}(t) =  H_x(\hat w(t),\hat p(t)),\q \bar H_u(\hat w(t),\hat p(t),\hat\lambda(t))=0 \q \mbox{a.e. in}\q [0,1].
$$
Notice that here and below, the dual variables $p$ and $\lambda$ are treated as row vectors, 
while $x$, $u$, $w$, $f$, and $G$ are treated as column vectors.


\subsection{ Second-order sufficient conditions for a weak local minimum }

Set 
$$ 
{ M}_{j}=\{t\in[0,1]: \; G_j(\hat u(t))=0\},\q j=1,\ldots,k.
$$  
Define the {\em critical cone}
$$\ba{rcl} K :=\Big\{\,w\in\W&:& 
\dot x(t)=f'(\hat w(t))w(t) \; \mbox{a.e. in}\; [0,1],\;  F'(\hat q )q \le0, \\  [2mm]
&& 
G'_j(\hat u(t))u(t)\le0 \;\;\mbox{a.e. on}\; M_j,\q   j=1,\ldots,k \,\Big\}.
\ea$$
It can be easily verified that  $ F'(\hat q )q =0$ for any element $w$ of the critical cone, and, moreover,
\be\label{cc} 
\ba{rcl} K =\Big\{\,w\in\W&:& 
\dot x(t)=f'(\hat w(t))w(t) \q \mbox{a.e. in}\q [0,1], \\ [2mm] &&   H_u(\hat w(t),\hat p(t))u(t)=0 \q 
\mbox{a.e. in}\q [0,1],\\ [2mm]
&& 
G'_j(\hat u(t))u(t)\le0 \;\;\mbox{a.e. on}\; M_j,\q   j=1,\ldots,k \,\Big\}.
\ea\ee
In what follows, we will use  representation (\ref{cc}) of the critical cone.

In many cases (say, in  problems of mathematical programming and calculus of variations) 
it is sufficient for local minimality that the critical cone consists only of the zero element. 
However, this is not the case for optimal control problems with a control constraint of the type $u(t)\in U$. 
Let us show this for the following problem which is somewhat different from (\ref{1})--(\ref{3}), since 
the dynamics is non-stationary. But formally, this problem can be reduced to  (\ref{1})--(\ref{3}) by introducing 
an additional state variable: $\dot y = t,\q y(0) = 0$. We leave the  verification of this fact to the reader. 

\begin{Example}
	Let $n=m=k=1$.  Consider the problem
	\bda
	& \min \{x(1) - x(0) \} &\\
	& \dot x(t) = tu(t)-(u(t))^2, \q u(t)\ge0 & \q \mbox{a.e. in}\q [0,1].
	\eda
	Set $\hat u=0$, $\hat x(t) = \hat x(0)$. The optimality system is satisfied with $\hat p(t) = 1$,
	$\hat \lambda(t)=t$. Here $K=\{0\}$. However,  $\hat u$ is not a  weak local minimum, because
	for the sequence
	$$
	u_s(t)=\left\{\ba{rcl}\frac{1}{s} &&
	0\le t\le \frac{1}{s}\vspace{2ex}\\ 0,&& \frac{1}{s} < t\le 1,
	\ea\right.
	$$ 
	we have $J(u_s) = - 1/(2s^3) < 0$ for all $s=1,2,\ldots$, and $\| u_s - \hat u \|_\infty \to 0$. 
Thus, the  condition $K=\{0\}$ is not sufficient for local minimality in this problem.
\end{Example}

	Let us give another equivalent representation of the critical cone. Set
$$ 
M^+(\hat\lambda_j)=\{t\in[0,1]: \; \hat\lambda_j(t)>0\}, \q j=1,\ldots,k.
$$ 
Then, due to (\ref{6}),
\be\label{ccc1}
\ba{rcl} K =\Big\{\,w\in\W &:& \dot x(t)=f'(\hat w(t))w(t)\q \mbox{a.e. in}\q [0,1], \\ [2mm] &&   G'_j(\hat u(t))u(t)\le0
\;\;\mbox{a.e. in}\;\; M_j; \\ [2mm] 
&& G'_j(\hat u(t))u(t)=0 \;\;\mbox{a.e. in}\;\;   M^+(\hat\lambda_j), \\ [2mm] &&  j=1,\ldots,k\,\Big\}.
\ea\ee 

We introduce an extension of the critical cone. For any $\delta>0$ and $j=1,\ldots,k$ we set
{$$ 
	M^+_\delta(\hat\lambda_j)=\{t\in[0,1]: \; \hat\lambda_j(t) >\delta\}.
	$$}
Next, for any $\delta>0$  we set
\be\label{cccd} \ba{rcl} 
K_\delta =\Big\{\,w\in\W&:&  \dot x(t)=f'(\hat w(t))w(t)\q \mbox{a.e. in}\q [0,1], \\ [2mm] &&  G'_j(\hat u(t))u(t)\le0
\;\;\mbox{a.e. on}\; M_j, \\ [2mm] &&  G'_j(\hat u(t))u(t)=0
\;\;\mbox{a.e. on}\;   M^+_\delta(\hat\lambda_j),\\ [2mm] &&  \;\;  j=1,\ldots,k \,\Big\}.
\ea\ee
Notice that the cones $K_\delta$ form a {non-increasing} family as $\delta\to0+$. {In particular, 
	$K \st K_\delta$ for any $\delta > 0$.}

Define the  {\em quadratic form}:
\be\label{12}
\Omega(w):= \langle F''(\hat q)q,q\rangle +
\int_0^1 \langle  \bar H_{ww}(\hat w(t),\hat p(t), \hat\lambda(t))w(t),w(t)\rangle\dd t,
\ee
{where} $ q=(x(0),x(1))$. 

\medskip       

\begin{Assumption}\label{assum4.3} There exist $\delta>0$ and $c_\delta>0$ such that 
	\be\label{13}
	\Omega(w)\ge c_\delta\big(|x(0)|^2 +\|u\|_2^2\big)\q  \forall\, w\in K_\delta.
	\ee
\end{Assumption}

\begin{Remark}
	Assumption \ref{assum4.3} is equivalent to the following: 
	there exist $\delta>0$ and $c_\delta>0$ such that 
	\be\label{14}
	\Omega(w)\ge c_\delta\big(\|x\|_\infty^2 +\|u\|_2^2\big)\q  \forall\, w\in K_\delta.
	\ee
\end{Remark}

Let us recall the following theorem, first published by the first author in the 1975 paper \cite{Osm-75}    
in a slightly different formulation.

\begin{Theorem}\label{th4.1}(sufficient second order condition) Let  Assumptions \ref{assum4.1}, 
	\ref{assum4.2}, \ref{assum4.3}  be fulfilled.
	Then there exist {$\e > 0$} and $c>0$ such that
	\be\label{gr} 
	J(w)-J(\hat w)\ge c\big(\|x-\hat x\|_\infty^2+\|u-\hat u\|^2_2 \big)
	\ee
	for all admissible $w=(x,u)\in W^{1,1}\times L^\infty$ such that  $\|w-\hat w\|_\infty<{\e}$.
\end{Theorem}

Proofs of the results in this subsection can be found, for example, in \cite{Osm-22}.

\subsection{An equivalent form of the second-order sufficient condition for local optimality}

Now we show that Assumption \ref{assum4.3} can be reformulated in terms of the critical cone $K$, 
instead of $K_\delta$, provided that an additional condition of Legendre type is fulfilled.

Let  $(\hat w,\hat p, \hat \lambda) \in \W \times  W^{1,1} \times L^\infty$,  and let  Assumptions
\ref{assum4.1}   and \ref{assum4.2} hold.

\begin{Assumption}\label{assum4.4}
	There exists  $c_0>0$  such that 
	\be\label{13a}
	\Omega(w)\ge c_0\big(|x(0)|^2 +\|u\|_2^2\big)\q  \forall\, w\in K.
	\nonumber
	\ee
\end{Assumption}

Further, for any $\delta>0$ and any $t\in[0,1]$ denote by $\C_\delta(t)$ the cone of all vectors $v\in\R^m$ 
satisfying for all 
$j=1,\ldots,k$ the conditions
$$  \C_\delta(t):=\q
\left \{ v\in\R^m\;:\q  \ba{l}  G_j'(\hat u(t))v\le 0 \q\mbox{if} \q G_j(\hat u(t))=0,\\ [2mm]
G_j'(\hat u(t))v=0 \q\mbox{if} \q \hat \lambda_j(t)>\delta, \q j=1,\ldots,k\ea 
\right\}. 
$$

For any $\delta>0$ and any $j\in \{1,\ldots,k\}$ we set
$$ 
m_\delta(\hat \lambda_j) :=\{t\in [0,1] :\; 0<\hat \lambda_j(t)\le\delta\}, \q    
m_\delta : =\bigcup_{j=1}^k  m_\delta(\hat \lambda_j).  
$$
Clearly, $\meas  m_\delta\to 0$ as $\delta\to0+$.

\begin{Assumption}\label{assum4.5}(strengthened Legendre condition on $m_\delta$). 
	There exist $\delta>0$ and $c_{\delta}^L>0$  such that for a.a. $t\in m_\delta$ we have 
	\be\label{13leg} 
	\langle \bar H_{uu}(\hat w(t),\hat p(t), \hat \lambda(t))  v,  v\rangle  \ge c_{\delta}^L  |v|^2\q  \forall\, v\in \C_\delta(t).
	\ee
\end{Assumption}

\begin{Proposition} Assumptions \ref{assum4.4} and \ref{assum4.5} together are equivalent 
	to Assumption \ref{assum4.3}.
\end{Proposition}

Thus, instead of Assumption \ref{assum4.3} we can use  Assumptions \ref{assum4.4} and \ref{assum4.5}  
in the sufficient second-order conditions of Theorem \ref{th4.1}.

The connection between the strengthened Legendre condition and the so-called  
``local quadratic growth of the Hamiltonian" (defined below) was studied by Bonnans and Osmolovskii 
in the 2012  paper \cite{B-Osm-2012}.

\begin{Definition}\label{def4.1} We say that the local quadratic growth condition of the Hamiltonian is fulfilled if 
	there exist $c_H>0$,  $\delta>0$, and $\varepsilon>0$ such that for  a.a. $t\in m_\delta  $ we have 
	$$ 
	H(\hat x(t),u,\hat p(t))-H(\hat x(t),\hat u(t),\hat p(t))\ge c_H|u-\hat u(t)|^2
	$$ 
	{for all} $u\in U$ {such that}  $|u-\hat u(t)| <\varepsilon $.
\end{Definition}

\begin{Proposition} Assumption \ref{assum4.5}  implies the local quadratic 
	growth condition of the Hamiltonian.
\end{Proposition}	

The converse is not true: the condition of the local quadratic growth 
of the Hamiltonian is somewhat finer than  Assumption \ref{assum4.5}. 

There is the following more subtle second-order sufficient condition for a weak local minimum at 
the point $\hat w$ 
in problem (\ref{1})-(\ref{3}).

\begin{Theorem} (sufficient second order condition) 
	Let the Assumptions \ref{assum4.1}, \ref{assum4.2},    \ref{assum4.4}   and  the local quadratic growth condition 
	of the Hamiltonian (Definition \ref{def4.1})  be satisfied. 	Then there exist {$\e_0 > 0$} 
and $c>0$ such that
	\be\label{gr1} 
	J(w)-J(\hat w)\ge c\big(\|x-\hat x\|_\infty^2+\|u-\hat u\|^2_2 \big) \nonumber
	\ee
	for all admissible $w=(x,u)\in W^{1,1}\times L^\infty$ such that  $\|w-\hat w\|_\infty<{\e_0}$. 
\end{Theorem}

Proofs of the results in this subsection can be found in \cite{Osm-22}.


\subsection{Strong metric subregularity of the optimality system} 

Now we formulate our main result concerning the SMSR property in optimal control: the optimality mapping 
associated with problem (\ref{1})-(\ref{3})  is strongly metrically subregular at a reference solution
$(\hat w,\hat p,\hat \lambda) = (\hat x,\hat u,\hat p,\hat \lambda) \in \W \times W^{1,1} \times L^\infty$ 
of the optimality system 
(\ref{3a})--(\ref{10}), provided that Assumptions \ref{assum4.1}, \ref{assum4.2}, and \ref{assum4.3}  hold.


\smallskip

Consider the {\em  perturbed system\,} of optimality conditions
(\ref{3a})--(\ref{10}):
\be\label{14a}  \lambda\ge0, \q \lambda  (G( u)-\eta)=0,\ee
\be \label{15}   (-p(0), p(1))=     
F'(q)+\nu,\ee
\be \label{16}   
\dot p + p \,f_x(w)= \pi,\ee
\be \label{17} 
pf_u(w)+\lambda G'(u)=\rho,
\ee
\begin{equation}
\label{20}	- 
\dot x + f(x,u) =\xi     
\end{equation} 
\be\label{21} 
G(u)\le \eta,\ee 
where  $\nu\in\R^{2n*}$,    $\pi\in L^1$,  
$\rho\in L^\infty$, $\xi\in L^1$,  $\eta\in L^\infty$. Note that $\nu$, $\pi$, 
and $\rho$ are treated as row vectors, while $\xi$ and $\eta$ are treated as column vectors.

Below we set
$$ 
\Delta x=x-\hat x,\q  \Delta u=u-\hat u,\q  \Delta w=(\Delta x,\Delta u)=
w-\hat w,\q \Delta p=p-\hat p, \q \Delta \lambda =\lambda-\hat \lambda,
$$ 
$$ 
\Delta q=(\Delta x(0),\Delta x(1))= (x(0)-\hat x(0), x(1)-\hat x(1))= (\Delta x_0, \Delta x_1),
$$ 
$$	\omega=(\nu,\pi, \rho,\xi, \eta),$$
$$  \|\omega\|:=|\nu|+ \|\pi\|_1+\|\rho\|_2+ \|\xi\|_1+\|\eta\|_2.
$$

\begin{Theorem}  Let  Assumptions \ref{assum4.1}, \ref{assum4.2}, and \ref{assum4.3}   be fulfilled.
	Then there exist reals {$\e > 0$}
	and $\kappa > 0$ such that if  
	\be\label{21p}  
	|\nu|+\|\pi\|_1+ \|\rho\|_\infty+	\|\xi\|_1 +	\|\eta\|_\infty \le {\e},
	\ee 
	then  for any solution $(x,u,p,\lambda)$ of the perturbed system (\ref{14a})--(\ref{21}) such that 
	$\|\Delta w\|_\infty\le\e$  the following estimates hold: 
	\be\label{21q} 
	\|\Delta x\|_{1,1}\le \kappa\|\omega\|,\q
	\|\Delta u\|_2 \le \kappa  \|\omega\|,\ee 
	\be\label{21r}  	\|\Delta p\|_{1,1} \le \kappa  \|\omega\|, \q \|\Delta\lambda \|_2\le \kappa \|\omega\|. 
	\ee
\end{Theorem}  

A proof of this theorem is given in \cite{Osm-Vel-23}.

\section*{Acknowledgments}

{This study is financed by the European Union-NextGenerationEU, through the
National Recovery and Resilience Plan of the Republic of Bulgaria, project
No
BG-RRP-2.004-0008-C01, and by the Austrian Science Foundation (FWF) under grant 
No I-4571-N.} 



\end{document}